\documentclass[10pt,twocolumn]{IEEEtran}
\usepackage{amssymb}  
\usepackage{amsfonts}
\usepackage{amsmath}

\usepackage{graphicx}      
\begin{document}

\title{Regions of Attraction for \\Hybrid Limit Cycles of Walking Robots \thanks{This work was supported by NSF Contract 0915148.}} 

\author{Ian R. Manchester, Mark M. Tobenkin, Michael Levashov, Russ Tedrake \\
CSAIL, Massachusetts Institute of Technology, USA.
\\\{irm, mmt, levashov, russt\}@mit.edu}

\maketitle
\begin{abstract}                
This paper illustrates the application of recent research in region-of-attraction analysis for nonlinear hybrid limit cycles. Three example systems are analyzed in detail: the van der Pol oscillator, the ``rimless wheel'', and the ``compass gait'', the latter two being simplified models of underactuated walking robots. The method used involves decomposition of the dynamics about the target cycle into tangential and transverse components, and a search for a Lyapunov function in the transverse dynamics using sum-of-squares analysis (semidefinite programming). Each example illuminates different aspects of the procedure, including optimization of transversal surfaces, the handling of impact maps, optimization of the Lyapunov function, and  orbitally-stabilizing control design.
\end{abstract}


\section{Introduction}

The purpose of this paper is to illustrate a new technique for estimation of regions of attraction for nonlinear hybrid limit cycles proposed in \cite{Manchester10a}. Three example systems have been chosen which illuminate different aspects of the method.

A major motivation for the work in this paper is control of underactuated ``dynamic walking'' robots (\cite{Collins05}). These robots can exhibit efficient, naturalistic, and highly dynamic gaits. However, control design and stability analysis for such robots is a challenging task since their dynamics are intrinsically hybrid and highly nonlinear. 

Estimates of regions of attraction can be useful for many problems, including analysis of different candidate control laws, generation of tree-based feedback motion-planning controllers (\cite{Tedrake10}), or planning transitions among a library of pre-stabilized walking gaits (\cite{Gregg10}).

The method involves choosing a decomposition of the dynamics into tangential and transversal components, and then searching for a Lyapunov function in the transversal components that verifies a ``tube'' about the limit cycle in which orbital stability is guaranteed. The verification is performed using sum-of-squares (SoS) programming.

The first example is the van der Pol oscillator, chosen because it is very well-known and well studied, and thus provides a good test of the method. With this example, we illustrate the importance of selecting the transversal decomposition intelligently. Indeed, it is shown that the commonly-used orthogonal transversal surfaces are often a bad choice.

The second is the ``rimless wheel'', a one-degree-of-freedom hybrid mechanical system, which serves as a simple model of a walking robot. The dynamics of the sytem are simple enough that much can be said about the rimless wheel analytically. Despite being simple, it exhibits hybrid (switching) behaviour representing the foot fall of a walking robot.

The third example is the ``compass-gait'' walker, a more complex model of a walking robot. For this system, one cannot derive analytical regions of attraction, so a computational approach is essential. This system requires transversal surface optimization, proper handling of impacts, and also orbitally-stabilizing control design based on a transverse linearization.

\subsection{Background}
The most well-known tool for analysing limit cycles is the Poincar\'e map: orbital stability is characterized by stability of an associated ``first-return map'', describing the repeated passes of the system through a single transversal hypersurface. Often a linearization of the first-return map is computed numerically, and its eigenvalues can be used to verify local orbital stability. Since the system's evolution is only analyzed on a single surface, regions of stability in the full state-space are difficult to evaluate via the Poincar\'e map.

A related technique known variably as ``transverse coordinates'' or ``moving Poincar\'e sections'' also has a long history and was certainly known to exist by Poincar\'e, however has not been much used in applications until recently due to difficulty in the relevant computations (see \cite{Hale80}). With this technique, a new coordinate system is defined on a family of transversal hypersurfaces which move about the orbit under study. In most cases, it is also used to study local stability, however it can be adapted to characterize regions of stability in the full state space.

The method we propose is to construct the transverse dynamics in regions of the orbit, and then  utilize the well-known sum-of-squares (SoS) relaxation of polynomial positivity which is amenable to efficient computation via semidefinite programming  (see, e.g., \cite{Shor87, Parrilo00, Parrilo03a}). The sum-of-squares relaxation has been previously used to characterize regions of stability of equilibria of nonlinear systems (see, e.g., 
\cite{Topcu08, Tan08, Henrion05}) and as a tool for constructive control design in \cite{Tedrake10}. 

There is comparatively little work on computing regions of stability of limit cycles. The proposed method has aspects in common with the surface Lyapunov functions proposed in \cite{Goncalves05}, however that method was restricted to piecewise linear systems.  The technique of cell-to-cell mapping, proposed by  \cite{Hsu80},  improves the efficiency of exhaustive grid-based methods of regional analysis and has been used in analysis of simple walking robots (\cite{Schwab01}), however the computational cost is still exponential in the dimension of the system.

\section{Problem Statement}

We consider the following class of hybrid systems with planar switching surfaces:
\begin{eqnarray}\label{eqn:sys_u1}
\dot x &=& f(x,u),  \ x\notin S^-\\
x^+&=&\Delta(x), \ x\in S^-.
\end{eqnarray}
where $x\in\mathbb R^n$, and $u\in\mathbb R^m$.
Suppose $f(\cdot)$ and $\Delta(\cdot)$ are smooth and $\Delta:S^-\rightarrow S^+$ where
\begin{eqnarray}
S^- &=& \{x: c_-'x=d_-\},\\
S^+ &=& \{x: c_+'x=d_+\},
\end{eqnarray}
$c_-, c_+\in \mathbb R^n$, and $d_-, d_+ \in \mathbb R$. Suppose $x^\star(\cdot)$ is a non-trivial $T$-periodic solution that undergoes $N$ impacts at times $\{t_1, t_2, ..., t_N\}+kT$ for integer $k$. 
We will assume that the impacts are not ``grazing'', i.e. $c_-'f(x^\star(t_i)) \ne 0$ and $c_+'f(x^\star(t_i)) \ne 0 $ for all $i$.

It is well-known that periodic solutions of autonomous systems cannot be asymptotically stable, since perturbations in phase are persistent. The more appropriate notion is {\em orbital stability} (see, e.g., \cite{Hale80, Hauser94, Shiriaev08}).

The analysis objective is to efficiently compute a region of state space $D\subset \mathbb R^n$ from which orbital stability to $x^\star(\cdot)$ is guaranteed.

\section{Transverse Dynamics and Regions of Orbital Stability}

In this section we briefly describe the process of verifying regions of orbital stability. Full details of each step are given in \cite{Manchester10a}.

The process we propose for finding regions of orbital stability is based on the construction of a smooth local change of coordinates $x \rightarrow (x_\perp, \tau)$. At each point $t\in [0, T]$ we define a hyperplane $S(t)$, with $S(0)=S(T)$, which is transversal to the solution $x^\star(\cdot)$, i.e. $\dot x^\star(t)\not\in S(t)$. 

The transversal surfaces are defined by:
\[
S(\tau)=\{y\in\mathbb R^n : z(\tau)'(y-x^\star(\tau))=0\}
\]
where $z:[0,T]\rightarrow\mathbb R^n$ is a vector function which will be optimized.

\begin{figure}[t]
\centering
\hskip -10pt\includegraphics[width=0.65\columnwidth]{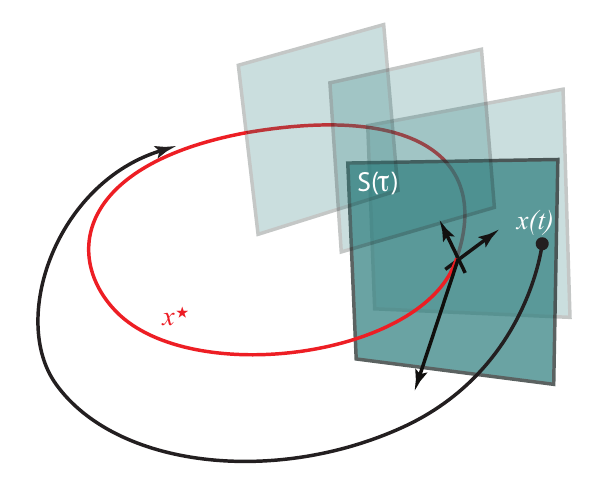}\\
\hskip 20pt\includegraphics[width=0.66\columnwidth]{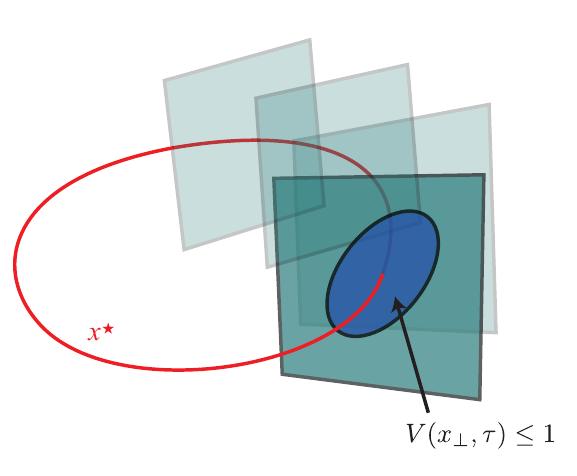}
\caption{Top: transversal surfaces $S(\tau)$ around the target orbit $x^\star$, with a particular solution $x(t)$ converging to the orbit $x^\star$. Bottom: a Lyapunov function defined on a transversal surface.}
\label{fig:ts1}
\end{figure}

Given a point $x$ nearby $x^\star(\cdot)$, the scalar $\tau \in [0, T)$ represents which of these transversal surfaces $S(\tau)$ the current state $x$ inhabits; the vector $x_\perp \in \mathbb R^{n-1}$ is the ``transversal'' state representing the location of $x$ within the hyperplane $S(\tau)$, with $x_\perp = 0$ implying that $x=x^\star(\tau)$. This is visualised in Figure \ref{fig:ts1}. 

The process for computing regions of orbital stability is as follows:
\begin{enumerate}
\item Select a family of transversal surfaces $S(\tau)$, and the associated transformation $x\rightarrow (x_\perp, \tau)$ such that at impact times the transversal surfaces line up with switching surfaces.
\item Compute the nonlinear dynamics in this new coordinate system as well as a periodic linear system representing the dynamics of $x_\perp$ close to the orbit: the {\em transverse linearization}.
\item Construct a candidate quadratic Lyapunov function associated with the transverse linearization via standard techniques from linear control theory.
\item Using this result as an initial seed, iteratively solve a sequence of sum-of-squares programs to compute maximal regions in which a Lyapunov function can be found verifying both well-posedness of the change of coordinates and orbital stability for the true nonlinear dynamics.
\end{enumerate}
For each $\tau\in[0,T]$, $\Pi(\tau)\in\mathbb R^{(n-1)\times n}$ is a smooth matrix function of $\tau$ projecting $x\rightarrow x_\perp$.

The dynamics in the new coordinates $x_\perp, \tau$ are given by the continuous dynamics
\begin{eqnarray}
\dot x_\perp&=& \dot \tau\left[\frac{d}{d\tau}\Pi(\tau)\right]\Pi(\tau)'x_\perp+\Pi(\tau)f(x^\star(\tau)+\Pi(\tau)'x_\perp)\notag\\ &&\,\,\,-\Pi(\tau)f(x^\star(\tau)) \dot \tau,  \label{eqn:xperpdot}\\
\label{eqn:taudot}
\dot \tau &=& \frac{z(\tau)'f(x^\star(\tau)+\Pi(\tau)'x_\perp)}{z(\tau)'f(x^\star(\tau))-\frac{d z(\tau)}{d \tau}'\Pi(\tau)'x_\perp}, 
\end{eqnarray}
for $t\ne t_i$, and impact updates
\begin{equation}\label{eqn:xperp_update}
x_\perp^+=\Pi(\tau_i^+)\big[\Delta_i\big(x^\star(\tau_i^-)+\Pi(\tau_i^-)'x_\perp^-\big)-x^\star(\tau_i^+)\big], 
\end{equation}
when $t=t_i$. The change of coordinates and the above dynamics are well-defined in a region around the target orbit $x^\star(\cdot)$.

\subsection{Regions of Orbital Stability}\label{sec:ver} 
Suppose there exists a Lyapunov function $V$ such that $V(x_\perp, \tau)>0, x_\perp \ne 0, V(0,\tau)=0$ for all $\tau\in [0,T]$ for which the following conditions hold on the level set $\{x:V(x_\perp, \tau)\le 1\}$:
\begin{eqnarray}
\frac{d}{dt}V(x_\perp, \tau)&\le&-\delta|x_\perp|^2,\label{eqn:hybridc1}\\
z(\tau)'f(x^\star(\tau))-\frac{\partial z(\tau)}{\partial \tau}'\Pi(\tau)'x_\perp&>&0,\label{eqn:hybridc2}
\end{eqnarray}
for some $\delta>0$.
 The first constraint verifies stability via the decreasing Lyapunov function; the second verifies well-posedness of the change of variables. If hybrid dynamics are considered, at the switching times one must verify the condition:
\begin{eqnarray}
V\Big(\Pi(\tau_i^+)\big[\Delta_i\big(x^\star(\tau_i^-)+ \Pi(\tau_i^-)'x_\perp\big) &&\notag\\ -x^\star(\tau_i^+)\big], \tau_i^+\Big) -V(x_\perp, \tau_i^-)&\le&0, \label{eqn:hybridd1}
\end{eqnarray}

If these conditions hold then the set $\{x:V(x_\perp, \tau)\le 1\}$ is an inner (conservative) estimate of the region of attraction to the limit cycle. Furthermore, for polynomial dynamics these can be relaxed to a SoS program. In some cases, one must also check that the impact surface is not reached before it is expected (see \cite{Manchester10a}).

Typically we sample a sufficiently fine finite sequence of $\tau\in[0, T]$ and verify the conditions on $x_\perp$ at each sample. In the optimization procedure, one must search for both $V(x_\perp, \tau)$ and Lagrange multipliers verifying the regional conditions. This problem is non-convex, however when fixing $V(x_\perp, \tau)$ and searching over Lagrange multipliers it is a semidefinite program, and when fixing Lagrange multipliers and searching over $V(x_\perp, \tau)$ it is a semidefinite program. Thus if one has a reasonable initial guess for $V(x_\perp, \tau)$ an iterative procedure can be applied. Whilst this is not guaranteed to find a global optimum, in practice the authors have found it works very well.

\subsection{Transverse Linearization}

In the construction of initial candidate Lyapunov functions and feedback controllers, we will make use of the linearization of the transverse dynamics:
\begin{eqnarray}\label{eqn:translin_u}
\dot x_\perp &=& A(t)x_\perp(t) +B(t)\bar u(t), \ \ t\ne t_i\\
x_\perp^+&=&A_dx_\perp, \ \ t = t_i.
\end{eqnarray}
representing the parts of \eqref{eqn:xperpdot} and \eqref{eqn:xperp_update} linear in $x_\perp$. Note that this can be given analytically for any system in the class.

\section{The van der Pol Oscillator}

The first example we consider is the van der Pol oscillator, defined by the following differential equation:
\[
\ddot y - \mu(1-y^2)\dot y+y=0.
\]
where $\mu$ is a constant, which we take equal to 1.
It is well-known that this system has a single unstable equilibrium at the origin, and a periodic cycle which is the limit from every other point in the plane. Thus, it makes a good simple example on which to test the proposed method. The above differential equation can obviously be rewritten in terms of a state $x =[x_1 \ x_2]':= [y \ \dot y]'$ in the form
\[
\dot x= f(x), \ f(x) = \begin{bmatrix}x_2\\(1-x_1^2)x_2-x_1\end{bmatrix}.
\]
To compute the periodic solution $x^\star$ we simply simulated the system forward from an initial condition away from the origin until it converged. Let $T$ be its period.

Given a vector $z(\tau)$, for a planar system the projection $\Pi(\tau)$ onto surfaces orthogonal to $z(\tau)$ is simple: $\Pi(\tau) = [-z_2(\tau) \ z_1(\tau)]$.
We start with transversal surfaces orthogonal to the system motion, i.e. $z(\tau) = f(x^\star(\tau))/|f(x^\star(\tau))|$. 

A natural candidate for a Lyapunov function is the solution of the Lyapunov differential equation for the transverse linearization:
\begin{equation}\label{eqn:lyap_eq}
\dot P(t)+A(t)'P(t)+P(t)A(t)+Q=0.
\end{equation}
For any $Q>0$ a unique periodic solution $P(t) = P(t+T)>0$ exists, and suggests a Lyapunov function of $x_\perp P(\tau) x_\perp$, as was suggested in \cite{Hauser94}.

\begin{figure}
\centering
\includegraphics[width=0.9\columnwidth]{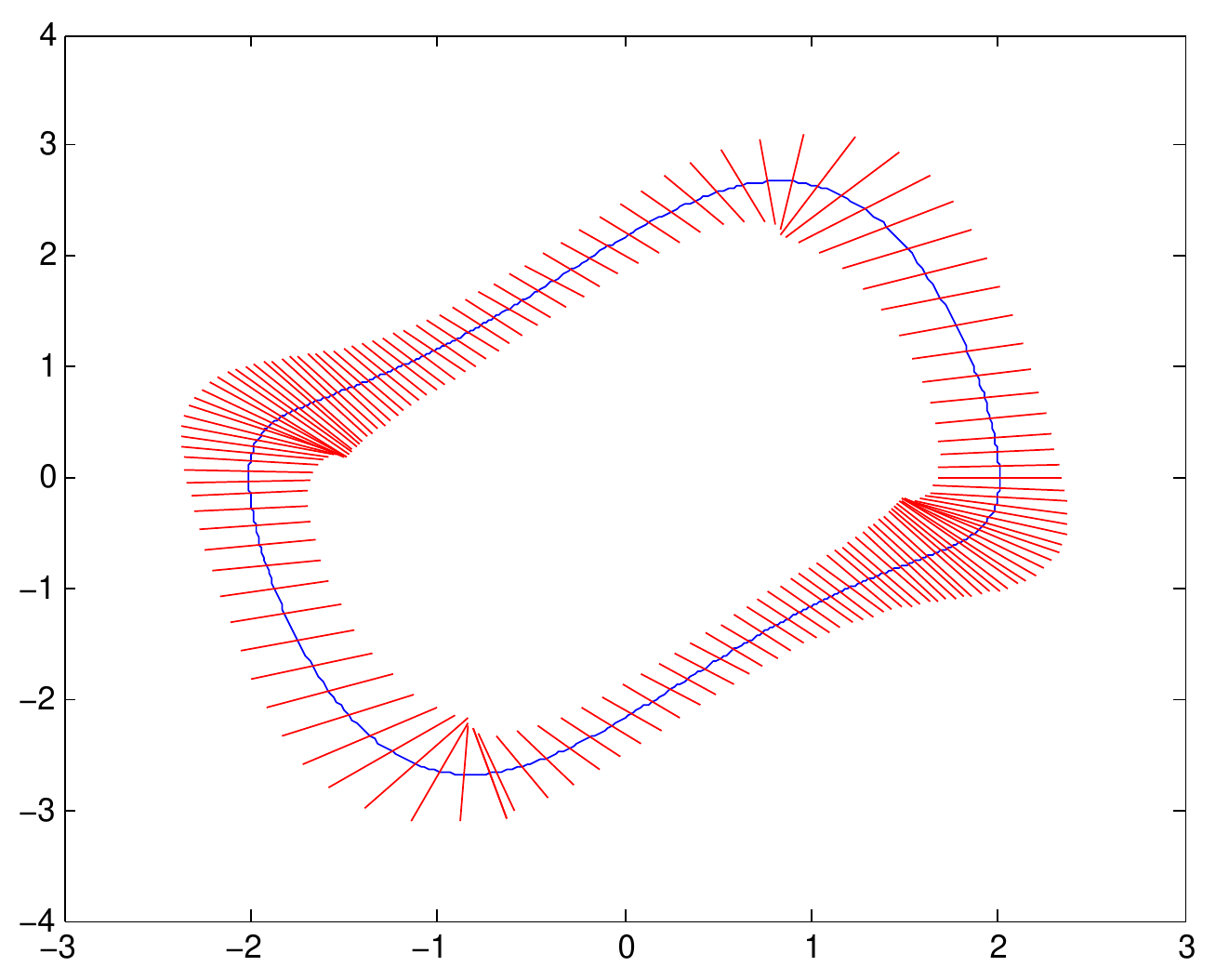}
\caption{Verified regions of orbital stability using constant rescalings of linear Lyapunov function  and orthogonal transversal surfaces.}
\label{fig:vdp1}
\end{figure}

We can then search for the maximal level set $x_\perp P(\tau) x_\perp \le \rho$ in which stability can be verified.  In the framework of Section \ref{sec:ver} this corresponds to verifying \eqref{eqn:hybridc1} and \eqref{eqn:hybridc2} for Lyapunov fuction of the form $V(x) = (1/\rho)x_\perp P(\tau) x_\perp$, which can be performed via a simple bisection search over $\rho$. The results are shown in Figure \ref{fig:vdp1}.

The regions are strongly limited by points at which the change of variables $x \rightarrow (x_\perp, \tau)$ becomes ill-defined. In the figure, these points resemble the hub of a bicycle wheel. Mathematically, they correspond to points at which the denominator of \eqref{eqn:taudot} goes to zero. This is clearly a consequence of the choice of transversal surfaces and motivates exploring other possible choices.

The van der Pol oscillator is a special case, since we know in advance that an orbital stability test must fail at the origin, the unstable equilibrium. With this in mind we construct radial transversal surfaces centred at the origin. Since $z(\tau)$ and $\Pi(\tau)$ are changed, we recompute the solution of \eqref{eqn:lyap_eq} and again search for the maximal level set of this Lyapunov function. 

\begin{figure}
\centering
\includegraphics[width=0.9\columnwidth]{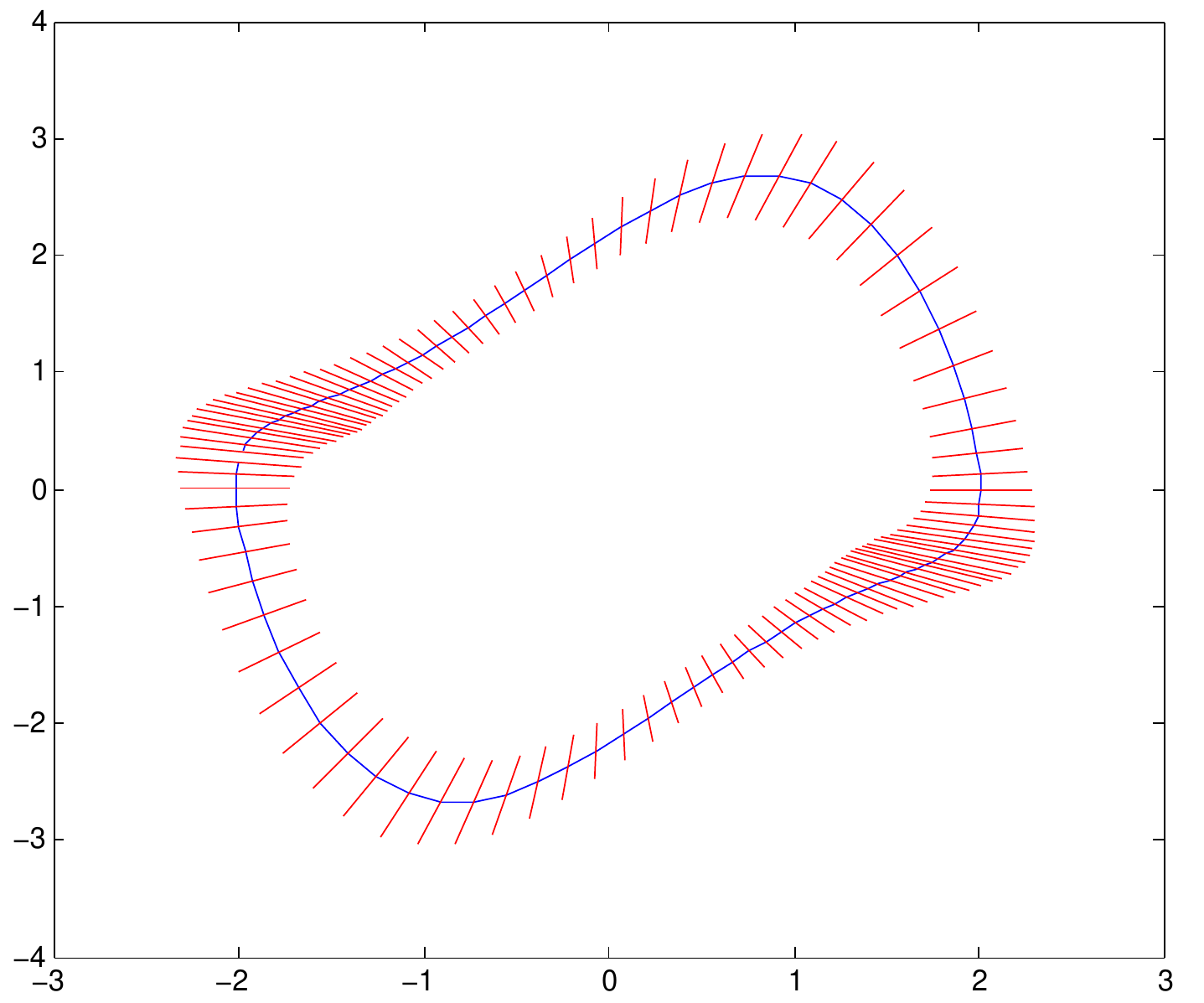}
\caption{Verified regions of orbital stability using constant rescalings of linear Lyapunov function and radial transversal surfaces.}
\label{fig:vdp2}
\end{figure}

The results are plotted in Figure \ref{fig:vdp2}. We no longer have the problem with singularities, however the regions are still quite thin. The reason is that whilst the Lyapunov ODE is guaranteed to {\em locally} verify orbital stability (see \cite{Hauser94}), it may not be a very good choice for verifying {\em regional stability}.

The result is much improved if we allow time-varying adjustment to the Lyapunov function via a scaling function $\sigma(\tau)>0$, i.e.
\begin{equation}\label{eqn:sig_rescale}
V(x_\perp,\tau) = \sigma(\tau)x_\perp'P(\tau)x_\perp
\end{equation}
where $\sigma(\cdot): [0, T] \rightarrow \mathbb R$, $\sigma(0)=\sigma(T)$, and $P(\tau)$ is as above.

\begin{figure}
\centering
 \includegraphics[width=0.9\columnwidth]{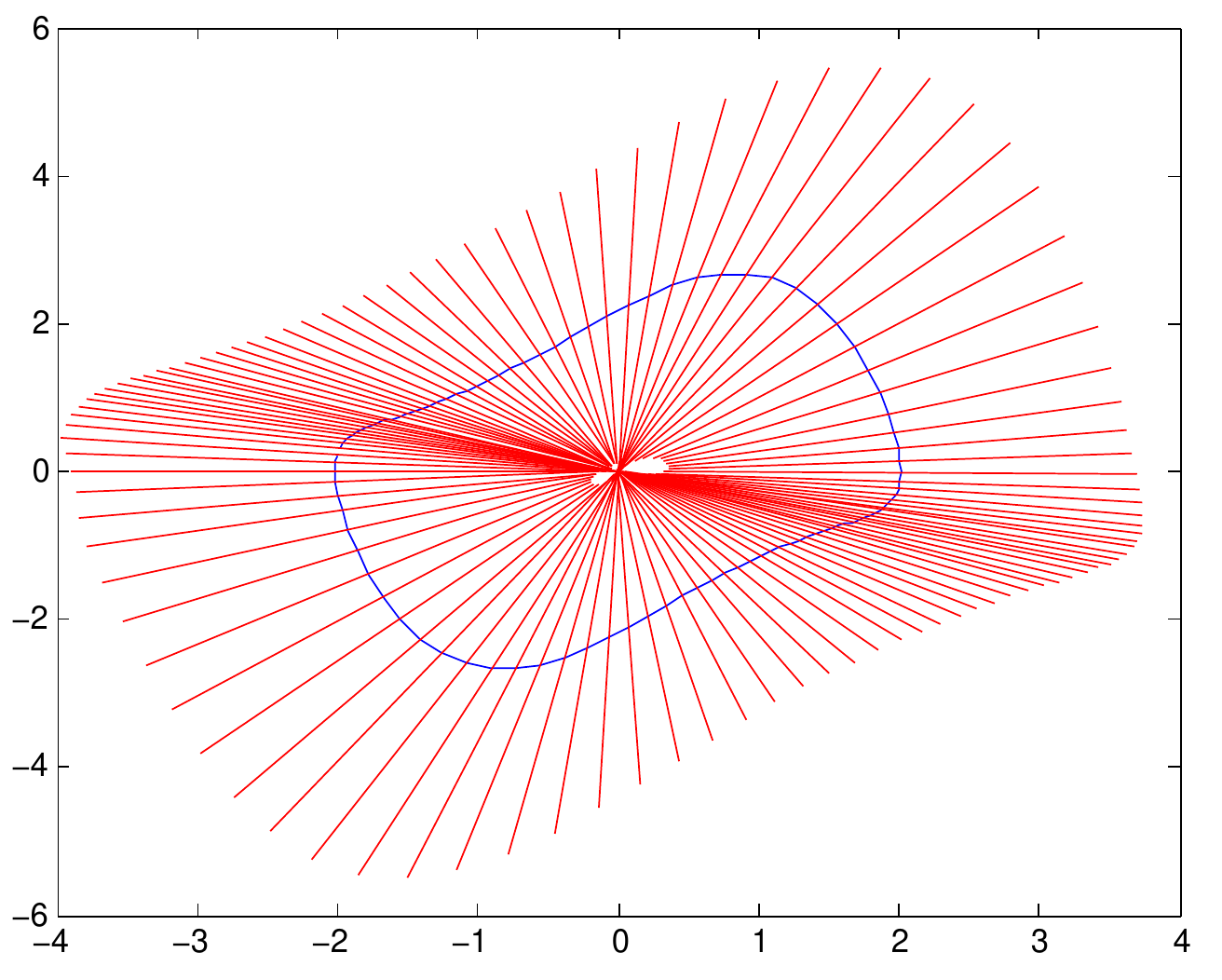}
\caption{Verified regions of orbital stability using time-varying $\sigma(\tau)$ and radial transversal surfaces.}
\label{fig:vdp3}
\end{figure}

The results are plotted in Figure \ref{fig:vdp3}. This figure is computed with $\sigma(\tau)$ a Bezier polynomial of order 20. It is observed that increasing the order of $\sigma(\tau)$ grows the region essentially to the origin. 

\begin{figure}
\centering
\includegraphics[width=0.9\columnwidth]{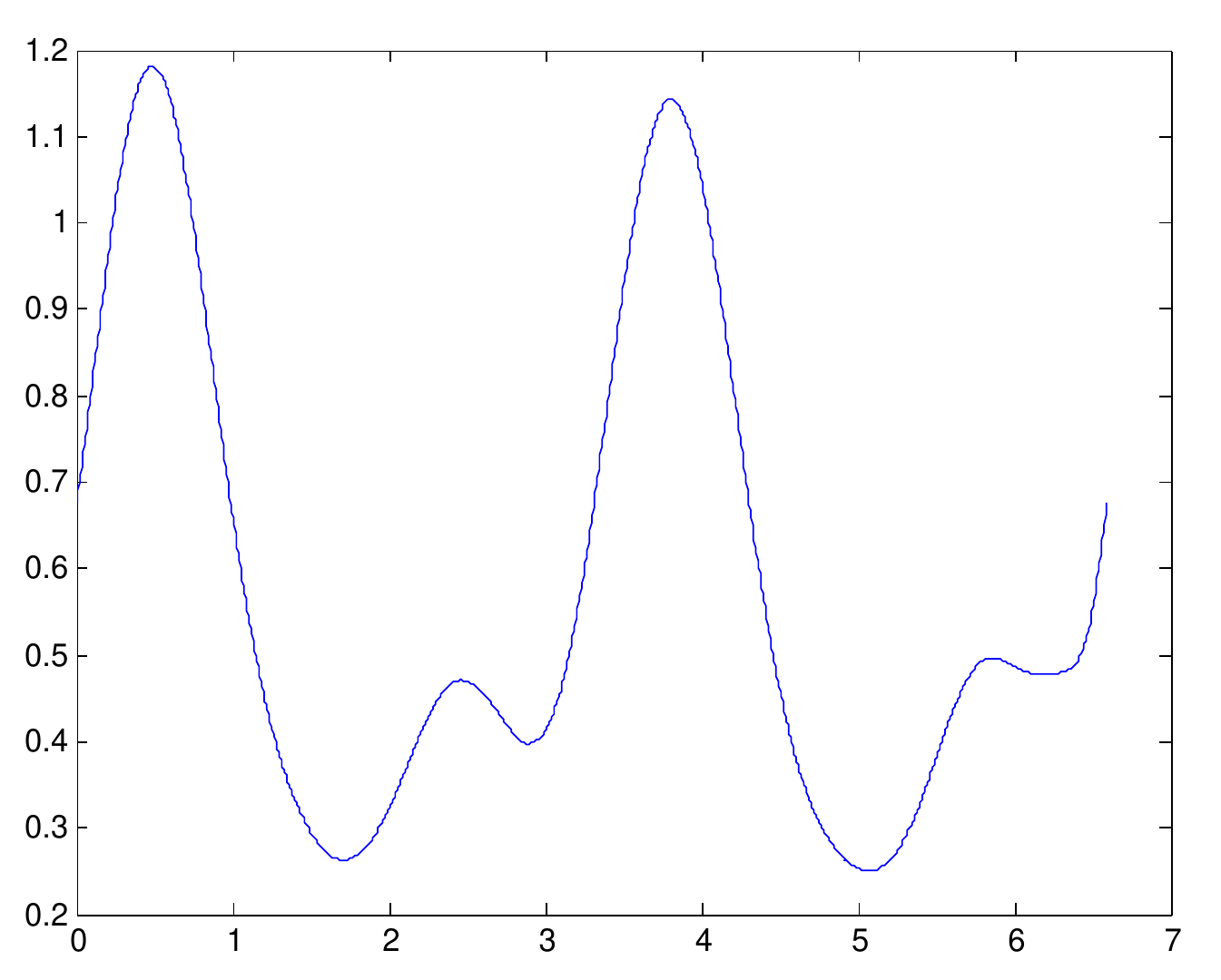}
\caption{The function $\sigma(\tau)$ used to verify the region in Fig  \ref{fig:vdp3}.}
\label{fig:sigma}
\end{figure}

Note that since we are searching over a class of candidate Lyapunov functions which are symmetric with respect to the orbit, the computed region is the best that can be achieved. The time-varying rescaling that achieved this is shown in Figure \ref{fig:sigma}.

\begin{figure}
\centering
\includegraphics[width=0.9\columnwidth]{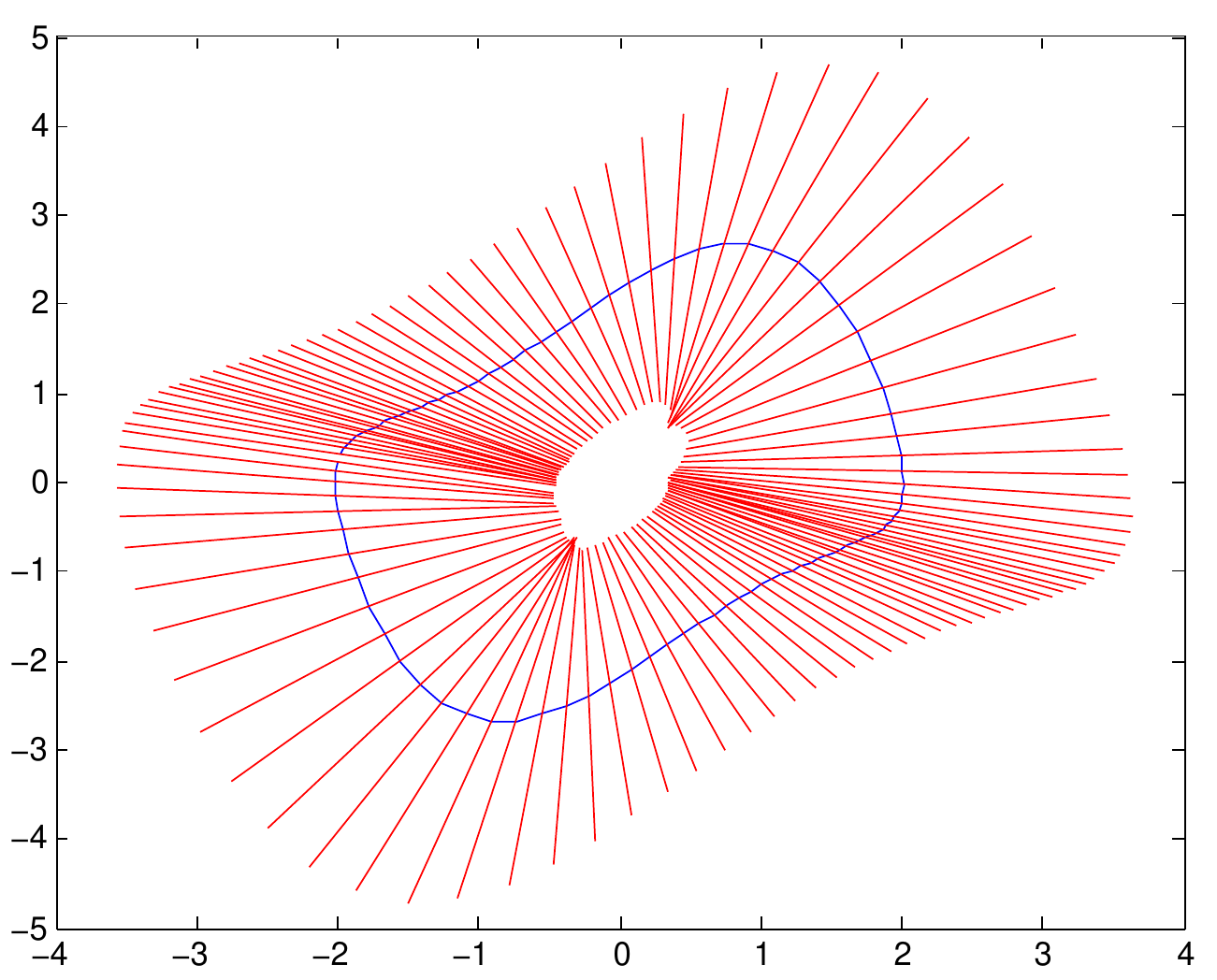}
\caption{Verified regions of orbital stability using time-varying $\sigma(\tau)$ and locally-optimized transversal surfaces.}
\label{fig:vdp_z}
\end{figure}

To obtain these results we made use of qualitative knowledge of the true region of stability in order to choose the transversal surfaces, which will not be possible for more complicated systems. In \cite{Manchester10a} a procedure for choosing $z(\tau)$ was suggested which depends only on local information about $f(x^\star, u^\star)$. The idea is to optimize the function  $z(\tau)$ so that the distance to the closest point of singularity is maximized. The results of applying this procedure are shown in Figure \ref{fig:vdp_z}.

This choice of transversal surfaces does not do quite as well the radial surfaces, which is not surprising. However it does substantially better than the orthogonal surfaces in Fig. \ref{fig:vdp1}. It seems there is still plenty of room for improvement in the optimization of $z(\tau)$, perhaps in an iterative procedure based on results of a prior region-of-stability computations.

\section{Rimless Wheel}

\begin{figure}
\centering
\includegraphics[width=0.6\columnwidth]{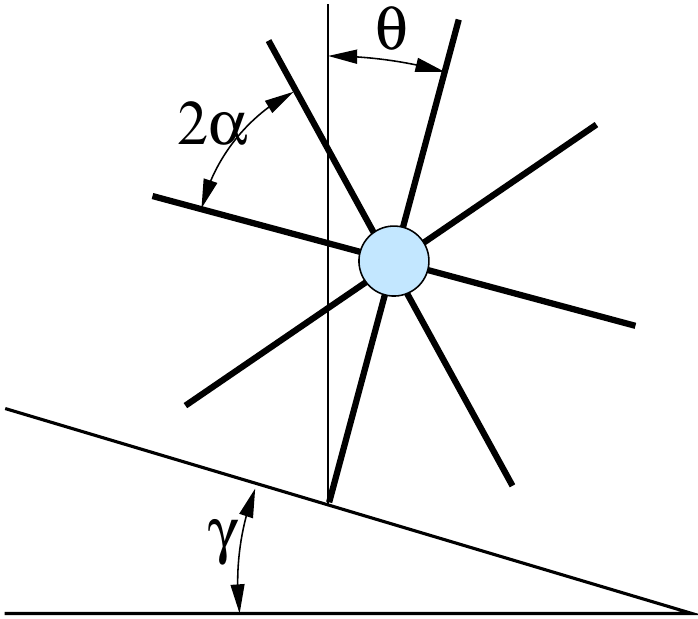}
\caption{Rimless wheel system}
\label{fig:rwheel_diag}
\end{figure}

The rimless wheel is a simple planar model of walking and consists of
a central mass with several `spokes' extending radially outward.
(See Figure \ref{fig:rwheel_diag}).
At any given moment one of the spokes is pinned at the ground,
and the system follows the dynamics of a simple pendulum,
$f(\theta,\dot\theta)=[\dot\theta, \sin(\theta)]'$.
When another spoke contacts the ground, the system undergoes an
inelastic collision governed by
$\dot\theta^+=\cos(2\alpha)\dot\theta$,
and the new spoke becomes the pinned one.

On a sufficiently inclined slope the system has a stable limit cycle,
for which the energy lost in collision is perfectly compensated by the
change in potential energy.
The rimless wheel has been analyzed in the literature and the basin
of attraction has been computed exactly (see \cite{Coleman98a}).
It is interesting to see how close our region of attraction estimation
method can approximate the actual basin for this system.

\begin{figure}
\centering
\includegraphics[width=1.1\columnwidth]{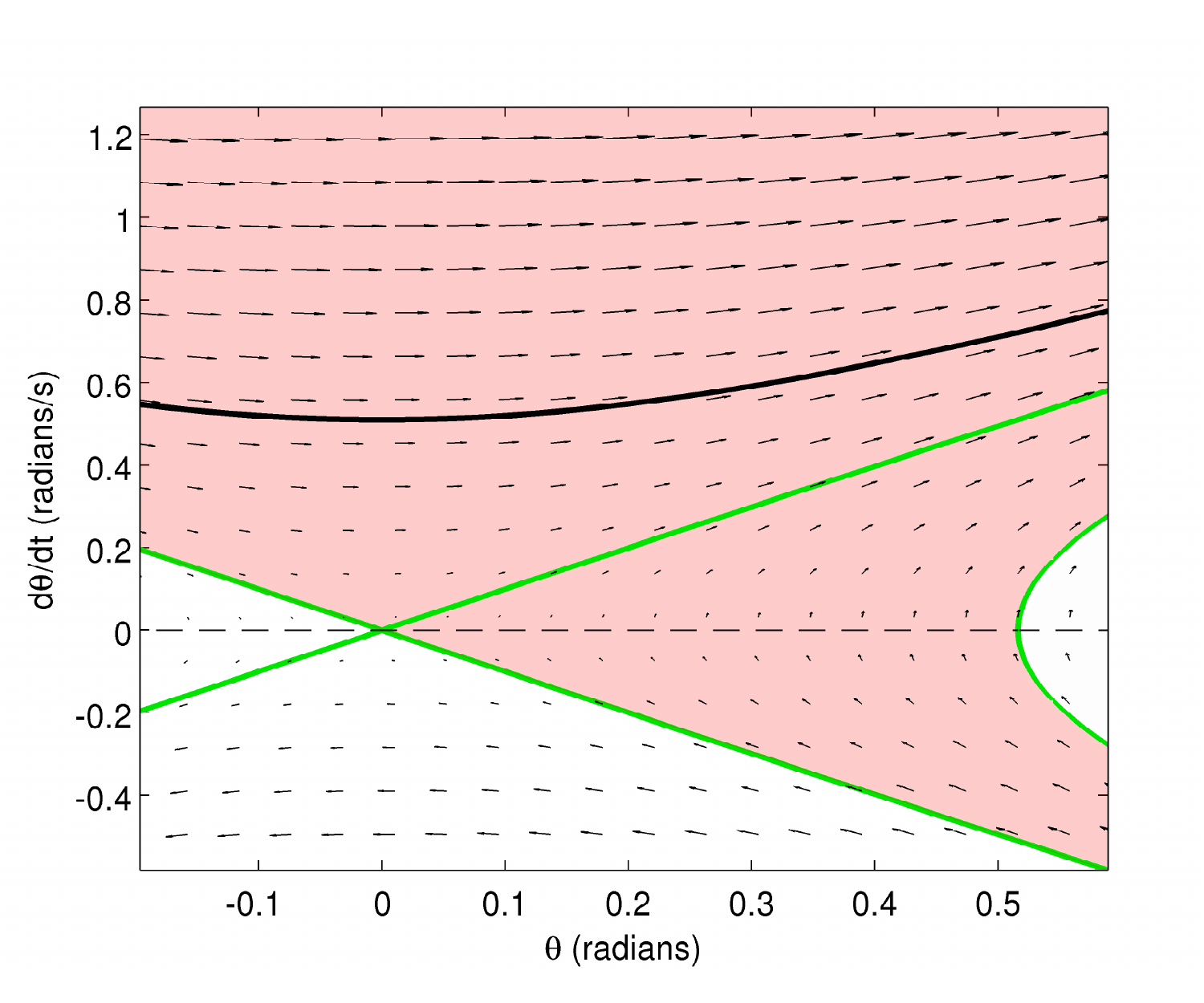}
\caption{Phase portrait of the rimless wheel system.
}
\label{fig:rwheel_basin}
\end{figure}

Figure \ref{fig:rwheel_basin} shows the phase portrait of the rimless wheel,
with arrows indicating the direction of the dynamics.
The right edge of the graph represents the collision surface that maps
to the left edge of the graph (or vice-versa, depending on the
direction of dynamics).
Because the impact depends only on the value of the angle, $\theta$, and
not on $\dot\theta$, the collision surfaces are vertical.
The thick green lines are the homoclinic orbits of the simple pendulum.
The thick black line shows the stable limit cycle,
and the shading shows a subset of its true region of attraction.

In this case, it is natural to select vertical transversal surfaces,
since there are no singularities in the change of variables, and the transversal surfaces are aligned with
switching surfaces.
The surfaces can be parametrized by $\theta$ as $\tau=\tau(\theta)$,
and the nominal trajectory is simply $\theta^\star(\theta)=\theta$
and $\dot\theta^\star(\theta)=\sqrt{2E-2\cos(\theta)}$, where $E$ is
the total energy of the system.
Then the transversal coordinate is the vertical position with respect to
the nominal trajectory: $x_\perp=\dot\theta-{\dot\theta}^\star(\theta)$.

With this choice of transversal surfaces, the dynamics in the new coordinate system are straightforward to compute:
\begin{eqnarray}
\dot{x}_\perp \notag
=\ddot\theta-\frac{d}{dt}\dot\theta^\star(\theta)
=\sin(\theta)- \frac{d \dot\theta^\star}{d \theta}(\theta)\dot\theta\\
=\sin(\theta)-\frac{\sin(\theta)}{\dot\theta^\star(\theta)}\dot\theta
=-\frac{\sin(\theta)}{\dot\theta^\star(\theta)}x_\perp.
\label{eqn:rw_xperpdot}
\end{eqnarray}
Interestingly, the transversal dynamics are linear for any given $\theta$.
Note that the transverse dynamics are also easily derived from equations
\eqref{eqn:xperpdot} and \eqref{eqn:taudot} by setting
$z(\tau)=[1, 0]'$ and $\Pi(\tau)=[0, 1]$.
This gives $\dot\tau=\frac{\dot\theta}{\dot\theta^\star(\theta)}$
and the equation \eqref{eqn:rw_xperpdot}, as expected.

As with the van der Pol oscillator, to find an initial candidate Lyapunov function we computed the unique periodic solution of the jump Lyapunov differential equation:
\begin{eqnarray}
-\dot P(t) &=& A(t)'P(t)+P(t)A(t) +Q, \ t \ne t_i \notag\\
P(\tau_i^-)&=&A_d(\tau_i)'P(\tau_i^+)A_d(\tau_i)+Q_i, \ t = t_i \notag
\end{eqnarray}
with $Q, Q_i>0$ and search over scalar rescalings:
$
V(x) = (1/\rho)x_\perp'P(\tau)x_\perp.
$

\begin{figure}
\centering
\includegraphics[width=1\columnwidth]{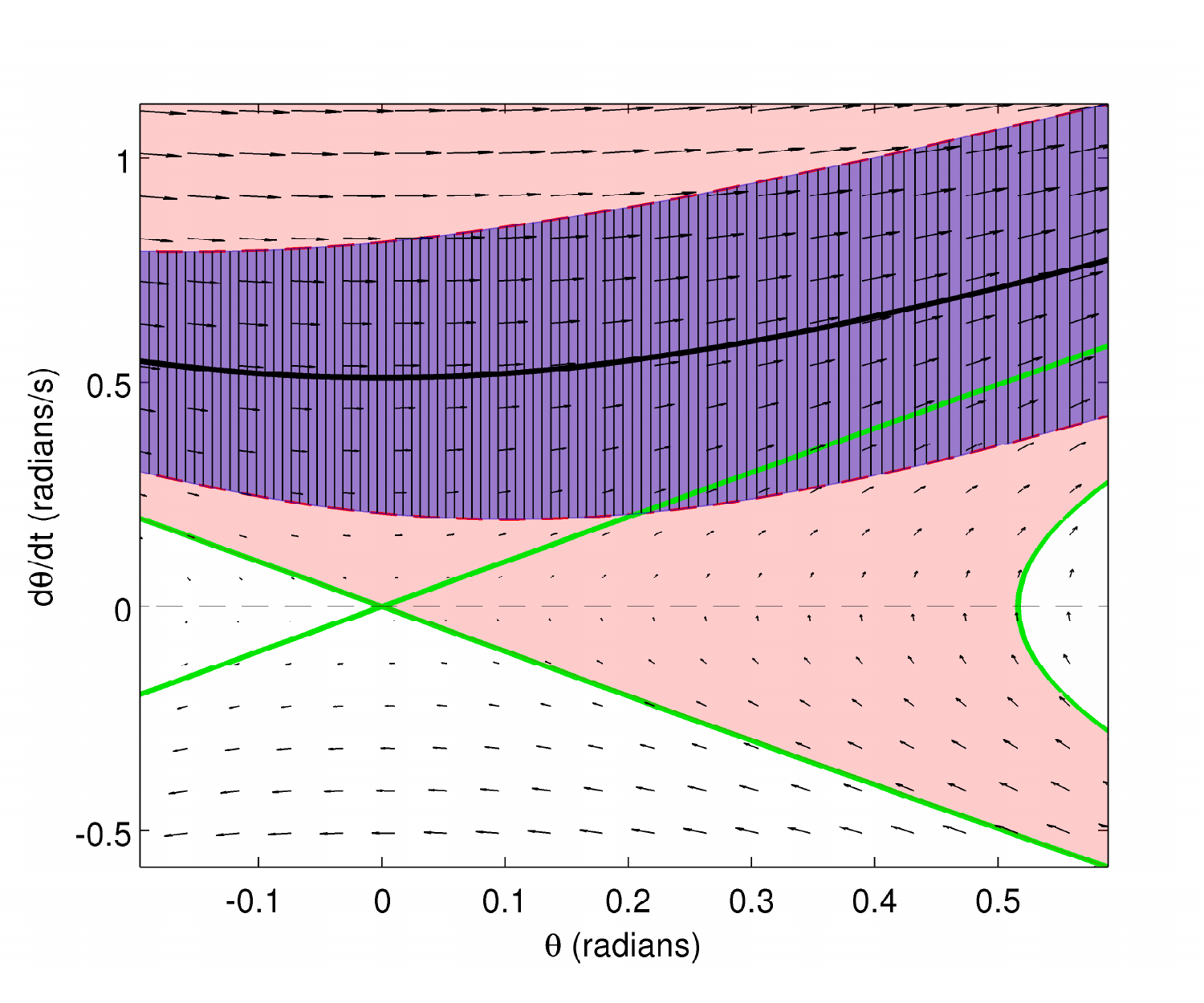}
\caption{Regions of Attraction for the rimless wheel limit cycle.
Light shaded region is the true RoA.
Dark shaded region is the verified RoA.
}

\label{fig:rwheel1}
\end{figure}

Figure \ref{fig:rwheel1} shows the results.
The discrete set of transverse surfaces can be seen as thin black
vertical lines around the orbit.
The computed basin of attraction (dark shading) is
within the true basin (light shading), but doesn't fill it fully.
This is not surprising, since we searched over a very restrictive set of
Lyapunov functions.

\begin{figure}
\centering
\includegraphics[width=1\columnwidth]{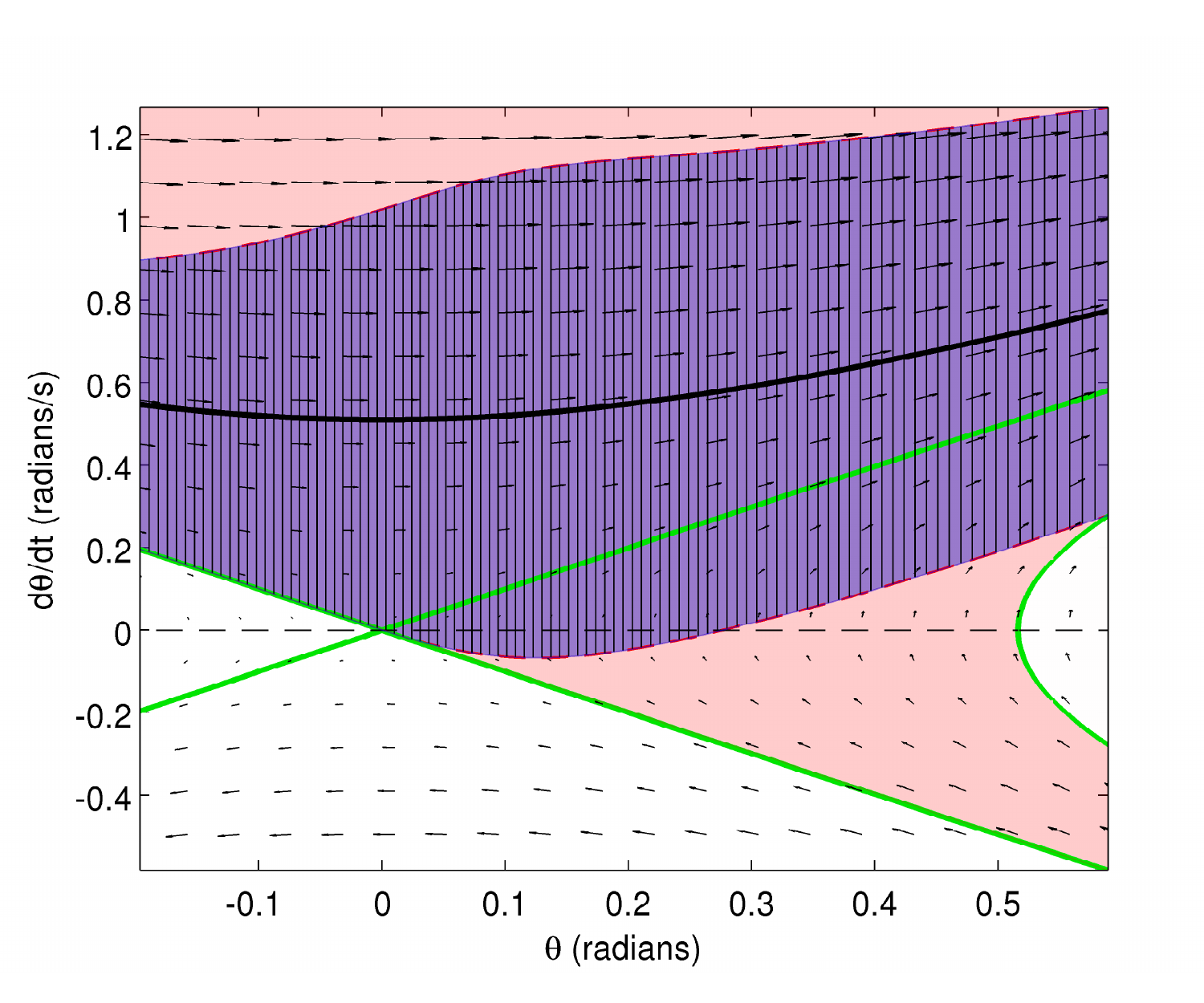}
\caption{Regions of Attraction for the rimless wheel limit cycle.
Light shaded region is the true RoA.
Dark shaded region is the verified RoA.
}

\label{fig:rwheel_verified}
\end{figure}

We then searched over $10^{\textrm{th}}$-order rescaling polynomials
$\sigma(\tau)$, as in \eqref{eqn:sig_rescale}, to maximize area of
the computed region.
Figure \ref{fig:rwheel_verified} shows the results and, as expected, the
basin of attraction is larger.

Note that the verified basin of attraction includes regions where
$\dot\theta < 0$ and, by the choice of transverse surfaces,
$\dot\tau < 0$.
The verification procedure still holds when the system state moves
backwards through the transverse surfaces.

\section{Compass-Gait Walker}

\begin{figure}
\centering
\includegraphics[width=0.6\columnwidth]{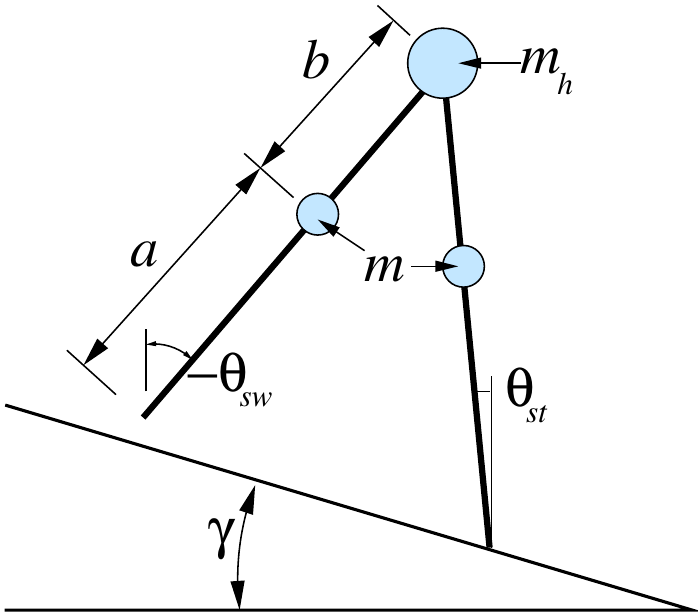}
\caption{Compass gait system}
\label{fig:cg_diag}
\end{figure}

The compass-gait walker is a two-degree-of-freedom (four-state)
nonlinear hybrid system shown in Figure \ref{fig:cg_diag}.
Similarly to the rimless wheel, one of the legs (the ``stance leg'') is always pinned
at the ground, and the other (the ``swing leg'') swings until it hits the ground, at which point the swing leg becomes the stance leg and vice versa.
A motor at the hip generates torque, $\tau$, between the legs.
Letting $q=[\theta_{sw},\theta_{st}]'$, $u=\tau$, and $l=a+b$
the continuous dynamics can be expressed in the standard manipulator form as
$$H(q)\ddot{q} + C(q,\dot{q})\dot{q} + G(q) = Bu,$$
where
\begin{gather*}
H = \begin{bmatrix}
mb^2~~ & ~~-mlb\cos(\theta_{st}-\theta_{sw}) \\
-mlb\cos(\theta_{st}-\theta_{sw})~~ & ~~(m_h+m)l^2 + ma^2
\end{bmatrix}, \\
C = \begin{bmatrix} 0~~ &  ~~mlb\sin(\theta_{st}-\theta_{sw})\dot\theta_{st} \\
-mlb\sin(\theta_{st}-\theta_{sw})\dot\theta_{sw}~~ & ~~0 \end{bmatrix}, \\
G = \begin{bmatrix} mbg\sin(\theta_{sw}) \\ -(m_hl + ma + ml)g\sin(\theta_{st}) \end{bmatrix},
 B = \begin{bmatrix} 1 \\ -1 \end{bmatrix}.
\end{gather*}
At the moments of impact, the coordinates undergo a relabling: $\theta_{sw}^+=\theta_{st}$ and $\theta_{st}^+=\theta_{sw}$.
The impact dynamics for velocities are derived assuming a perfect inelastic collision:
\[
Q^+_\alpha\dot q^+ = Q^-_\alpha\dot q^-
\]
where
\begin{eqnarray}
Q^-_\alpha &=& 
\begin{bmatrix}-mab & -mab+(m_hl^2+2mal)\cos(2\alpha)\\0&-mab\end{bmatrix}\notag
\end{eqnarray}
\begin{align}
&Q^+_\alpha= \notag\\
&\begin{bmatrix}mb(b-l\cos(2\alpha)) & ml(l-b\cos(2\alpha))+ma^2+m_hl^2\\mb^2&-mbl\cos(2\alpha)\end{bmatrix}\notag
\end{align}
and $\alpha = \frac{\theta_{sw}-\theta_{st}}{2}$. The impact takes place when $\theta_{sw}+\theta_{st}+2\gamma=0$, where $\gamma$ is the slope of the ground. Note that although the switching dynamics are highly nonlinear, the switching {\em surfaces} are planar, matching the assumptions of \cite{Manchester10a}.

For some combination of parameters and initial conditions the
system exhibits a limit cycle behavior of walking downhill passively
(with zero torque).
For this reason, the model is often used for studying bipedal walking
(see, e.g. \cite{Goswami98,Spong05, Westervelt07b,Byl08f, Manchester09} and many others).
Prior to this work, studies of the small basin of
attraction of this system were limited to exhaustive simulation.


We orbitally stabilize the limit-cycle using LQR control in the transverse coordinates, and analyze the region of attraction for the closed loop system.
We begin by computing the $T$-periodic stable limit cycle solution numerically from suitable initial conditions.  We next optimize for a set of transversal surfaces at $N=40$ sample points along the trajectory.  These surfaces are constrained to align with planar switching surfaces capturing the ``swing leg'' losing then regaining contact with the ground.  We construct the transverse coordinates, and the resulting transverse dynamics.

A set of transversal surfaces were chosen via the optimization procedure suggested in \cite{Manchester10a}.  An LQR controller was then designed for the transverse linearization by solving for the periodic positive-definite solution of the jump-Ricatti equation:
\begin{eqnarray}
-\dot P &=& A'P+SA -PBR^{-1}B'P +Q, \ t \ne t_i \notag\\
P(\tau_i^-)&=&A_d(\tau_i)'P(\tau_i^+)A_d(\tau_i)+Q_i, \ t = t_i \notag
\end{eqnarray}
with $R, Q, Q_i > 0$. The feedback control is given by:
\[
u(\tau, x_\perp) = -R^{-1}B(\tau)'P(\tau)x_\perp.
\]
Note that there is no nominal control command, since we are stabilizing a passive walking cycle.

To enable the use of SoS optimization, we approximate the non-polynomial dynamics of the walker via a third order Taylor expansion.  This expansion is performed around each of the $N$ sample points. 

\begin{figure}
\centering
\includegraphics[width=1\columnwidth]{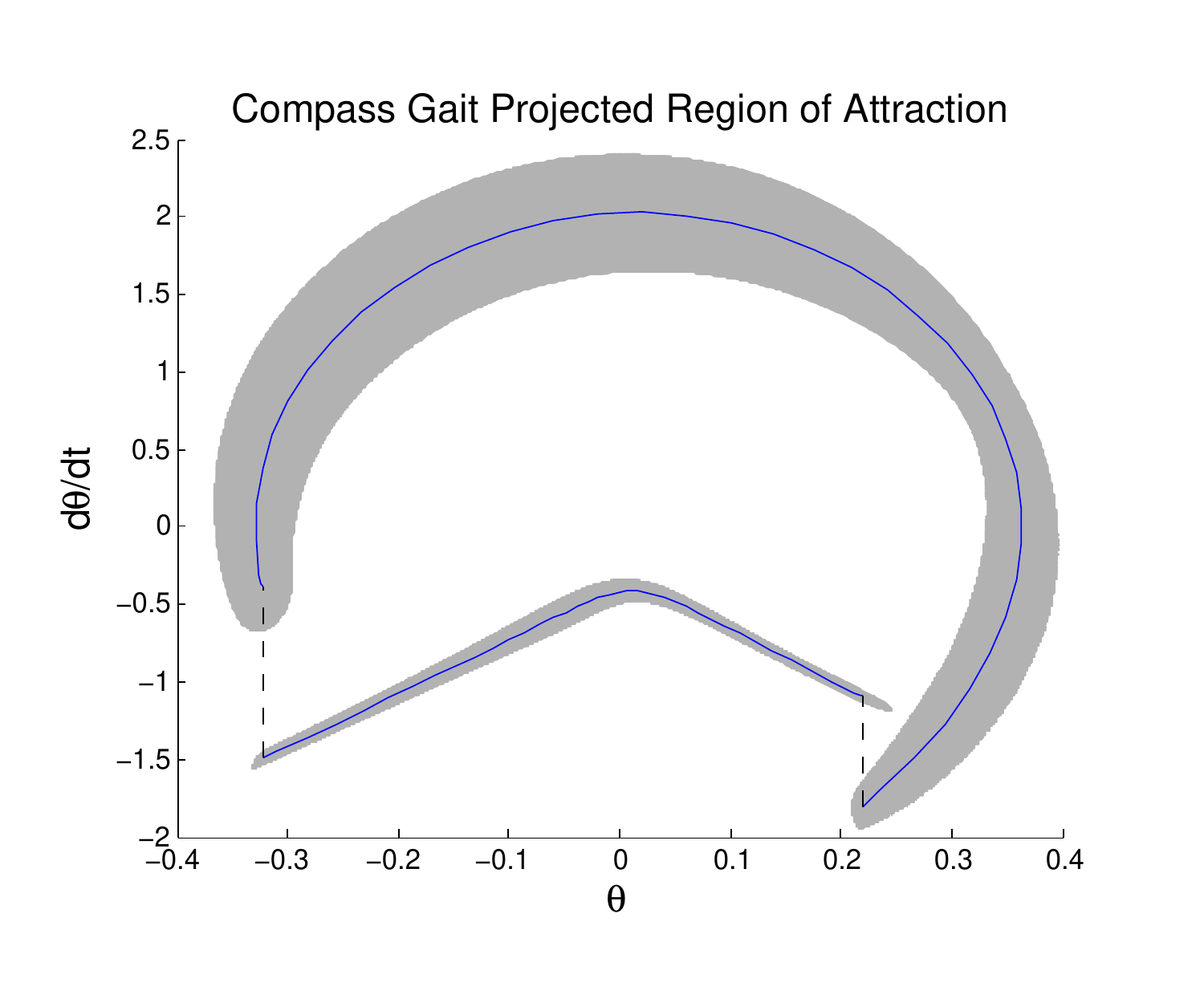}
\caption{Regions of Attraction for the compass-gait walker limit cycle.  We present two of six possible projections of the RoA.  The upper curve plots the trajectory and RoA projected into the $(\theta_1, \dot \theta_1)$ plane, that of the ``swing leg''.  The lower curve is the equivalent for the ``stance'' leg coordinates $(\theta_2, \dot \theta_2)$.}
\label{fig:cgw}
\end{figure}

The Riccati equation provides us with an initial candidate Lyapunov function:
$$V_0(x_\perp,\tau) = \sigma_0 x_\perp'P(\tau)x_\perp.$$
Taking a suitably small $\sigma_0$, we find a valid region of attraction.
To optimize the region of attraction, we rescale $V_0$ by a polynomial $\rho: [0,T] \mapsto (0,\infty)$:
$$V(x_\perp,\tau) = \frac{1}{\rho(\tau)}V_0(x_\perp,\tau).$$
We maximize the integral of $\rho(t)$ over the interval as a surrogate for the volume of the region of attraction.
For the compass gait system this parameterization proved better numerically conditioned than directly scaling the Lyapunov function.
The following observations allow our constraints to remain linear in the parameters of $\rho(t)$.  Note that:
\begin{align}
  \frac{d}{dt}V(x_\perp,\tau) = \frac{1}{\rho(\tau)^2}\Bigg[
  \rho(\tau)\frac{d}{dt} V_0(x_\perp,\tau) \nonumber
  -\frac{\partial \rho}{\partial \tau} \dot \tau V_0(x_\perp,\tau)
\Bigg],
\end{align}
To enforce
$\frac{d}{dt} V(x_\perp,\tau) \leq -\delta |x_\perp|^2$, we can exploit the fact that $\delta$ can be an arbitrary positive constant, and instead require:
\begin{flalign}
 \rho(\tau)\frac{d}{dt} V_0(x_\perp,\tau) -
  \frac{\partial}{\partial \tau} \rho(\tau) \dot \tau V_0(x_\perp,\tau) \leq \delta_0 |x_\perp|^2.\nonumber
\end{flalign}
As $\rho(t)$ is bounded above on $[0,T]$ an appropriate $\delta$ exists.
Finally, $\{ x \; | \; V(x_\perp,\tau) < 1\}$ is identically $\{ x \; | \; V_0(x_\perp,\tau) < \rho(\tau)\}$.
Figure~\ref{fig:cgw} presents overlapped plots of the region of attraction projected into the $(\theta_1, \dot \theta_1)$ and $(\theta_2, \dot \theta_2)$ planes.

These regions indicate that the controller can stabilize much larger variations in the swing-leg than in the stance-leg, which matches intuition on the compass gait walker.

\subsection{Accuracy of Taylor Expansion}

\begin{figure}
\centering
\includegraphics[width=1\columnwidth]{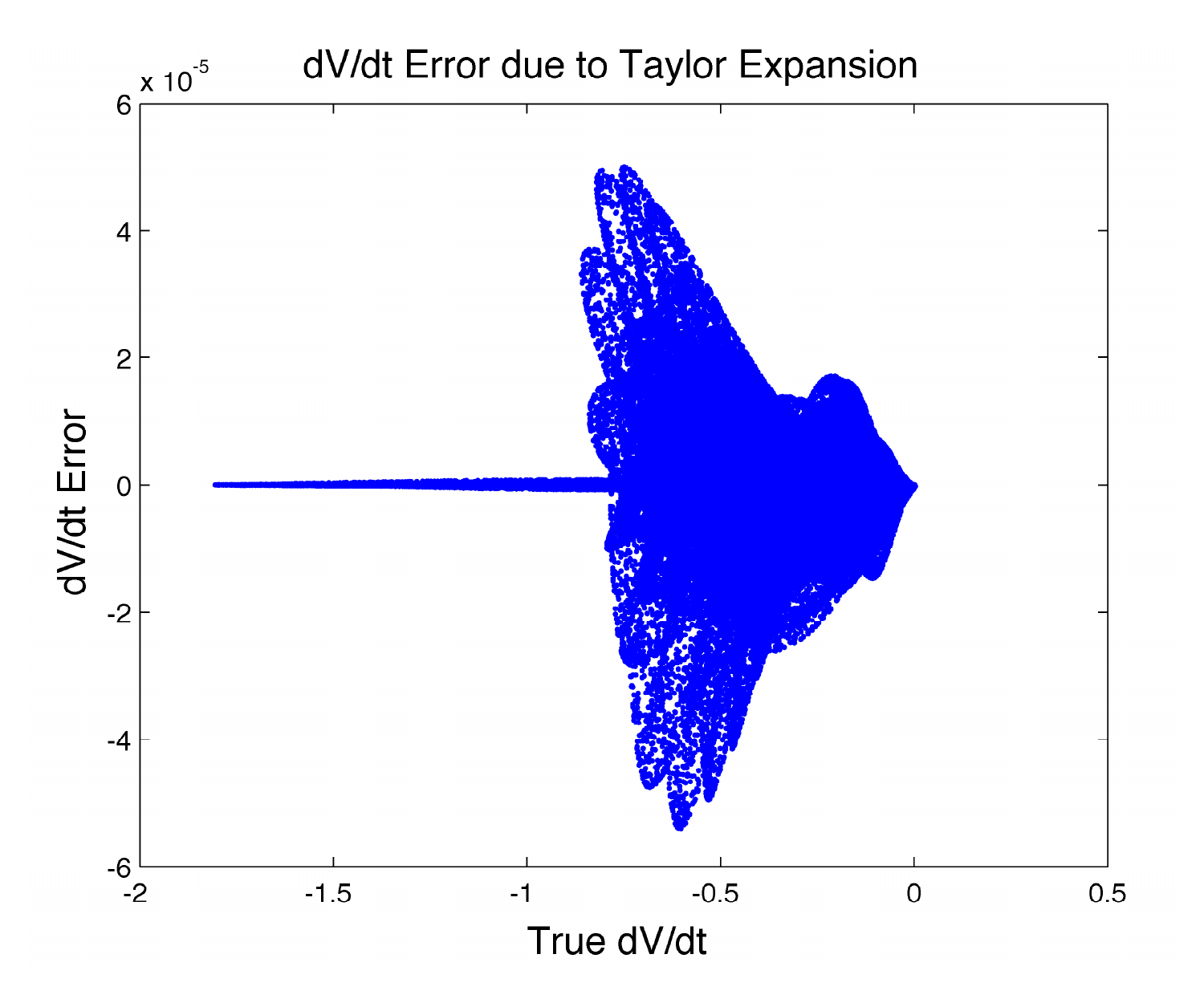}
\caption{The above compares
  the rate of change of the Lyapunov function for the true and Taylor expanded dynamics at the boundary of the approximate region-of-attraction. 99.7\% of samples were verified correctly.}
\label{fig:cgw-tay}
\end{figure}

Since we are approximating the dynamics by a Taylor expansion, it is important to check whether we can trust the verification.  To examine this approximation, we sample points on the boundary of the certified region of attraction.  For each point, we compute $\frac{d}{dt} V(x_\perp,\tau)$ with the true dynamics and Taylor expansion approximation.  At each time step we sampled 10,000 points from the boundary.  In Figure~\ref{fig:cgw-tay} we plot the true value of $\frac{d}{dt} V$ against the approximation error.  Note that the vertical scale of the plot is approximately $10^{-5}$ vs. an almost unit scaling along the horizontal, and essentially all the points are to the left of zero.

Of the 400,000 samples, 104 had $\frac{d}{dt} V > 0$, so the verification using Taylor series was correct at 99.97\% of samples. The largest value of  $\frac{d}{dt} V$ was 0.0026. Note that isolated samples of positive $\frac{d}{dt} V$ do not necessarily imply instability from that point, it is possible (indeed, likely) that the Lyapunov function was not appropriate for proving stability at that point.

\subsection{Computation Times}

In computing the regions of attraction, we alternate between optimizing over $V(x_\perp, \tau)$ -- a ``$V$-step'', and optimizing over Lagrange multipliers -- an ``$L$-step''. It took four $V$-$L$ iterations until covergence. The computations were performed on a 2.66 GHz intel Core i7 processor with 8GB of RAM, and timings for each iteration are shown in Table \ref{tab1}. The total computation time was around 17.7 minutes.

Note, however, that due to the sampling in $\tau$,  each $L$-step is made up of 41 independent optimizations (40 continuous samples and one impact map) each of which took between 1.6 and 2.6 seconds. Since these computations are trivially parallelizable, substantial speed increases could be achieved on multiprocessor machines.

\begin{table}\label{tab1}
\centering
\begin{tabular}{|c|c|c|c|c|}
\hline
Iteration& 1&2&3&4\\ \hline\hline
$L$-step time (s) &75.0 &90.2 &80.8&78.7 \\ \hline
$V$-step time (s) &194 &183 &201 & 170\\ \hline
\end{tabular}
\vskip5pt
\caption{Computation times for the compass-gait walker.}
\end{table}



\section{Summary and Future Work}

The purpose of this paper has been to demonstrate the application of a new method for estimating regions of attraction for nonlinear hybrid limit cycles. Three simple examples were chosen that elucidated important aspects of the technique, including selection of transversal surfaces, the handling of impacts, optimization of the Lyapunov function, and control design.

This work can be extended in a number of directions. It will serve as an essential component in recent algorithms for feedback control and motion planning in \cite{Tedrake10}. Implementation on more complex models and experimental validation will be an important test. 

In this paper we searched for regions via $\tau$-varying rescalings of quadratic Lyapunov functions from linear theory. This was done primarily to simplify the optimization procedure, however in principle it is possible to search over any polynomial Lyapunov functions, which would presumably give substantially bigger (and possibly asymmetric) regions for many systems.

There are many connections between sum-of-squares verification and robust control and stability theory via integral quadratic constraints (\cite{Megretski97, Petersen00}) which will be investigated.

\bibliographystyle{IEEEtran}
\bibliography{elib}

\end{document}